\documentclass[12pt]{article}

\usepackage[latin1]{inputenc}
\usepackage{bbm}
\usepackage{stmaryrd}
\usepackage{latexsym}
\usepackage{amsfonts}
\usepackage{amsmath,amssymb,amscd}
\usepackage{yfonts}
\usepackage{dsfont}
\usepackage{mathabx}
\usepackage{color}
\usepackage[dvips]{graphicx}
\usepackage{epsfig}
\usepackage{psfrag}

\DeclareFontFamily{U}{txsyc}{}
\DeclareFontShape{U}{txsyc}{m}{n}{
   <-> txsyc%
}{}
\DeclareFontShape{U}{txsyc}{bx}{n}{
   <-> txbsyc%
}{}
\DeclareFontShape{U}{txsyc}{l}{n}{<->ssub * txsyc/m/n}{}
\DeclareFontShape{U}{txsyc}{b}{n}{<->ssub * txsyc/bx/n}{}
\DeclareSymbolFont{symbolsC}{U}{txsyc}{m}{n}
\SetSymbolFont{symbolsC}{bold}{U}{txsyc}{bx}{n}
\DeclareFontSubstitution{U}{txsyc}{m}{n}
\DeclareMathSymbol{\df}{\mathrel}{symbolsC}{"42}
\DeclareMathSymbol{\fd}{\mathrel}{symbolsC}{"43}
\DeclareMathSymbol{\lJoin}{\mathrel}{symbolsC}{"58}
\DeclareMathSymbol{\rJoin}{\mathrel}{symbolsC}{"59}

\newcommand{\cC}{{\cal C}}
\newcommand{\cD}{{\cal D}}

\newcommand{\cL}{{\cal L}}
\newcommand{\cN}{{\cal N}}
\newcommand{\cM}{{\cal M}}
\newcommand{\cO}{{\cal O}}

\newcommand{\cT}{{\cal T}}

\newcommand{\LL}{\mathbb{L}}
\newcommand{\NN}{\mathbb{N}}

\newcommand{\RR}{\mathbb{R}}

\newcommand{\TT}{\mathbb{T}}

\newcommand{\di}{\displaystyle}
\newcommand{\iy}{\infty}
\newcommand{\lt}{\left}
\newcommand{\me}{\medskip}
\newcommand{\na}{\nabla}
\newcommand{\pa}{\partial}
\newcommand{\ri}{\rightarrow}
\newcommand{\rt}{\right}
\newcommand{\sm}{\smallskip}
\newcommand{\tr}{\triangle}

\newcommand{\wi}{\widetilde}

\newcommand{\fo}{\forall\ }
\newcommand{\Hess}{\mathrm{Hess}}
\newcommand{\lan}{\lt\langle}
\newcommand{\lve}{\lt\vert}
\newcommand{\lVe}{\lt\Vert}
\newcommand{\osc}{\mathrm{osc}}

\newcommand{\ran}{\rt\rangle}
\newcommand{\rve}{\rt\vert}
\newcommand{\rVe}{\rt\Vert}

\newcommand{\st}{\,:\,}

\newcommand{\un}{\mathds{1}}
\newcommand{\Var}{\mathrm{Var}}

\newcommand{\bq}{\begin{eqnarray*}}
\newcommand{\bqn}[1]{\begin{eqnarray}\label{#1}}
\newcommand{\eq}{\end{eqnarray*}}
\newcommand{\eqn}{\end{eqnarray}}
\newcommand{\wwtbp}{\par\hfill $\blacksquare$\par\me\noindent}
\newcommand{\thistitlepagestyle}{}

\newcommand{\ttsim}{\raise.17ex\hbox{$\scriptstyle\mathtt{\sim}$}}

\newtheorem{pro}{Proposition} 

\newtheorem{lem}[pro]{Lemma}
\newtheorem{theo}[pro]{Theorem}
\renewcommand{\thepro}{\arabic{pro}}

\newenvironment{rem}
{\par\me\refstepcounter{pro}\noindent{\bf Remark \thepro\ }}
{\par\hfill $\square$\par\me\noindent}
\newcommand{\proof}{\par\me\noindent\textbf{Proof}\par\sm\noindent}

\newcommand{\Na}{N^{(\alpha)}}

\title{A stochastic algorithm finding generalized means on compact manifolds}
\author{\Large Marc Arnaudon${}^\dagger$ and Laurent Miclo${}^\ddagger$
}
\date{\box1
 \box2
}

\begin{document}

\setbox1=\vbox{
\large
\begin{center}
${}^\dagger$ 
Institut de Mathématiques de Bordeaux, UMR 5251\\
Université de Bordeaux and CNRS, France\\
\end{center}
}
\setbox2=\vbox{
\large
\begin{center}
 ${}^\ddagger$
Institut de Mathématiques de Toulouse, UMR 5219\\
Universit\'e de Toulouse and CNRS, France\\
\end{center}
} 
\setbox3=\vbox{
\hbox{${}^\dagger$ marc.arnaudon@math.u-bordeaux1.fr\\}
\vskip1mm
\hbox{Institut de Mathématiques de Bordeaux \\}
\hbox{Université de Bordeaux 1\\}
\hbox{351, Cours de la Libération\\}  
\hbox{F33405 Talence Cedex, France
\\}
}
\setbox4=\vbox{
\hbox{${}^\ddagger$ miclo@math.univ-toulouse.fr\\}
\vskip1mm
\hbox{Institut de Mathématiques de Toulouse\\}
\hbox{Université Paul Sabatier\\}
\hbox{118, route de Narbonne\\} 
\hbox{31062 Toulouse Cedex 9, France\\}
}
\setbox5=\vbox{
\box3
 \vskip5mm
 \box4
}

\maketitle
\thistitlepagestyle
\abstract{
A stochastic algorithm is proposed,  finding the set of generalized means associated to a probability measure $\nu$ on a compact Riemannian manifold and a continuous cost function $\kappa$ on $M\times M$. Generalized means include $p$-means for $p\in(0,\infty)$, computed with any continuous distance function, not necessarily the Riemannian distance. They also include means for lengths computed from Finsler metrics, or for divergences.  

The algorithm is fed sequentially with independent random variables $(Y_n)_{n\in\NN}$ distributed according to $\nu$ and this is the only knowledge of $\nu$
required. It  evolves like a Brownian motion between the times it jumps
 in direction of the $Y_n$. Its principle is based on simulated annealing and homogenization, so that temperature and approximations schemes must be tuned up.  The proof relies on the investigation of the evolution of a time-inhomogeneous $\LL^2$ functional and on the corresponding spectral gap estimates due to Holley, Kusuoka and Stroock.
}
\vfill\null
{\small
\textbf{Keywords: }
Stochastic algorithms, simulated annealing, homogenization, probability measures on compact Riemannian manifolds, intrinsic means, instantaneous invariant measures, Gibbs measures,
spectral gap at small temperature.
\par
\vskip.3cm
\textbf{MSC2010:} first: 60D05
, secondary: 58C35, 60J75, 37A30, 47G20, 53C21, 60J65.
}\par

\newpage

\section{Introduction}

The purpose of this paper is to present a stochastic algorithm 
 finding the  generalized means of a probability measure $\nu$ defined on a compact manifold $M$. A Riemannian metric is used, only for the algorithm.
 
Algorithms for finding means, medians or minimax centers have been the object of many investigations, see e.g. \cite{Le:04}, \cite{Sturm:05}, \cite{Groisser:05}, \cite{Groisser:06}, \cite{Badoiu-Clarkson:03}, \cite{Yang:10}, \cite{Bonnabel:11}, \cite{Afsari-Tron-Vidal:11}, \cite{Arnaudon-al:12}, \cite{Cardot-Cenac-Zitt:12}, \cite{Arnaudon-Nielsen:12}. In these references a gradient descent algorithm is used, or a stochastic version of it avoiding to compute the gradient of the functional to minimize. Either the functional to minimize has only one local minimum which is also global, or (\cite{Bonnabel:11}) a local minimum is seeked. The case of Karcher means in the circle is investigated in  \cite{Charlier:11} and \cite{2011arXiv1108.2141H}. In this special situation the global minimum of the functional can be found by more or less explicit formula. 

For generalized means on compact manifolds the situation is different since the functional~\eqref{U} to minimize may have many local minima, and no explicit formula for a global minimum can be expected. In~\cite{Arnaudon-Miclo:13} the case of compact symmetric spaces has been investigated and a continuous inhomogeneous diffusion process has been constructed which converges in probability to the set of $p$-means. In~\cite{Arnaudon} the case of $p$-means on the circle is treated. A Markov process is constructed  which has Brownian continuous part and more and more frequent jumps in the direction of independent random variables with law $\nu$. It is proven that it converges in probability to the set of $p$-means of~$\nu$. Both~\cite{Arnaudon-Miclo:13} and~\cite{Arnaudon} use simulated annealing techniques.

The purpose of this paper is to extend the construction in~\cite{Arnaudon} to all compact manifolds, and to all generalized means. 

\par\me
So let be given $\nu$ a probability measure on $M$, a compact Riemannian manifold.
Denote by $\kappa : M\times M\to \RR$ a continuous function  and consider the continuous mapping
\bqn{U}
U\st M\ni x\mapsto \int_M \kappa(x,y)\, \nu(dy).\eqn
A global minimum of $U$ is called a $\kappa$-mean or generalized mean of $\nu$ and let $\cM$ be their set, which is non-empty in
the above compact setting.\\
In practice the knowledge
of $\nu$ is often given by a sequence $Y\df(Y_n)_{n\in \NN}$ of independent random variables, identically distributed
according to $\nu$.
So let us present a stochastic algorithm using this data and enabling to find some
elements of $\cM$. It is based on simulated annealing and homogenization procedures.
Thus we will need respectively an inverse temperature evolution $\beta\st \RR_+\ri \RR_+$, an inverse speed up evolution
$\alpha\st \RR_+\ri \RR_+^*$ and a regularization function $\delta\st \RR_+\ri \RR_+^*$. Typically, $\beta_t$ is non-decreasing, $\alpha_t$ and $\delta_t$ are non-increasing and we have $\lim_{t\ri+\iy}\beta_t=+\iy$,
$\lim_{t\ri+\iy}\alpha_t=0$ and $\lim_{t\ri+\iy}\delta_t=0$, but we are looking for more precise conditions so that the
stochastic algorithm we describe below finds $\cM$.
\\
For  finding $\cM$, we will use a regularization of $\kappa$ with the heat kernel $p(\delta,x,z)$. So define  for $\delta>0$   $\displaystyle \kappa_\delta(x,y)=\int_Mp(\delta,x,z)\kappa(z,y)\,\lambda(dz)$ where $\lambda$ denotes the Lebesgue measure (namely the unnormalized Riemannian measure) and 
\bqn{Us}
U_\delta\st M\ni x\mapsto \int_M \kappa_\delta(x,y)\, \nu(dy).\eqn

Notice any other regularization which would satisfy the estimates~\eqref{gradlnp} and~\eqref{partialnp} below and which would be easier to compute  could be used instead of the heat kernel.

Let $N\df(N_t)_{t\geq 0}$ be a standard Poisson process: it starts at 0 at time 0 and has jumps of length 1 whose interarrival times are independent and distributed according to exponential random variables of parameter 1.
The process $N$ is assumed to be independent from the chain $Y$.
We define the speeded-up process $\Na\df(\Na_t)_{t\geq 0}$ via
\bqn{Na}
\fo t\geq0,\qquad \Na_t&\df& N_{\int_0^t \frac1{\alpha_s}\,ds}.
\eqn
Consider the time-inhomogeneous Markov process $X\df(X_t)_{t\geq0}$ which evolves in $M$  in the following way:
if $T>0$ is a jump time of $\Na$, then $X$ jumps at the same time, from $X_{T-}$
to $X_T\df\phi_\delta(\beta_T\alpha_T,X_{T_-},Y_{\Na_{T}})$, where 
\bqn{phi}
s\mapsto\phi_\delta(s,x,y)\in M
\eqn
is the value at time $s$ of the flow started at $x$ of the vector field $z\mapsto -\frac12\na_z\kappa_\delta(\cdot,y)$. In particular 
\begin{equation}
\label{phiprime}
\phi_\delta'(s,x,y)=-\frac12\na\kappa_\delta(\cdot,y)(\phi_{\delta}(s,x,y))
\end{equation}
where $\phi_\delta'$ denotes the derivative with respect to the first variable.

 Typically we will have $\di\lim_{t\ri +\iy}\alpha_t\beta_t\|\na\kappa_{\delta_t}\|_\infty=0$,
so that for sufficiently large jump-times $T$, $X_T$ will be "between" $X_{T-}$ and $Y_{\Na_{T}}$ and quite close to $X_{T-}$. To proceed with the construction, we require that between consecutive jump times (and between time 0 and
the first jump time), $X$ evolves as a Brownian motion, relatively to the Riemannian structure of
$M$ (see for instance the book of Ikeda and Watanabe \cite{MR1011252}) and independently 
of $Y$ and $N$.
Informally, the evolution of the algorithm $X$ can be summarized by the It\^o equation (in centers of exponential charts)
\bq
\fo t\geq0,\qquad
dX_t&=&\sigma(X_t)d B_t+ \overrightarrow{X_{t-}\phi_{\delta_t}(\beta_t\alpha_t,X_{t_-},Y_{\Na_{t}})}\,d\Na_t\eq
where $(B_t)_{t\geq0}$ is a Brownian motion on some  $\RR^m$, for all $x\in M$ $\sigma(x) :\RR^m\to T_xM$ is linear satisfying $\sigma\sigma^\ast={\rm id}$, $x\mapsto \sigma(x)$ is smooth,  and where $(Y_{\Na_t})_{t\geq0}$ should be interpreted as a fast auxiliary process.
The law of $X$ is then entirely determined by the initial distribution $m_0=\cL(X_0)$.
More generally at any time $t\geq0$, denote by $m_t$ the law of $X_t$.\par
We will prove that the above algorithm $X$ finds  in probability at large times 
the set of means $\cM$.
Let us define a constant $b\geq 0$, coming  from the theory of simulated annealing (cf.\ for instance
Holley, Kusuoka and Stroock \cite{MR995752}) in the following way.
For any $x,y\in M$, let $\cC_{x,y}$ be the set of continuous paths $p\df (p(t))_{0\leq t\leq1}$
going from $p(0)=x$ to $p(1)=y$.
The elevation $U(p)$ of such a path $p$ relatively to $U$ is defined by
\bq
U(p)&\df& \max_{t\in[0,1]} U(p(t))\eq
and the minimal elevation $U(x,y)$ between $x$ and $y$ is
given by
\bq
U(x,y)&\df& \min_{p\in\cC_{x,y}} U(p).\eq
Then we consider
\bqn{bU}
b(U)&\df& \max_{x,y\in M}U(x,y)-U(x)-U(y)+\min_M U\eqn
This constant can also be seen as the largest depth of a well
not encountering
a fixed global minimum of $U$. Namely, if $x_0\in \cM$, 
then we have
\bqn{bU2}
b(U)&=& \max_{y\in M} U(x_0,y)-U(y)\eqn
independently of the choice of 
$x_0\in\cM$.

\par\sm
With these notations, the main result of this paper is:
\begin{theo}\label{t1}
For any scheme of the form
\bqn{ab}
\fo t\geq0,\qquad \lt\{\begin{array}{rcl}
\alpha_t&\df& (1+t)^{-1}\\
\beta_t&\df& c^{-1}\ln(1+t)\\
\delta_t&\df& \ln(2+t)^{-1}
\end{array}\rt.\eqn
where  $c>b(U)$, we have for any neighborhood $\cN$ of $\cM$ and for any $m_0$,
\bqn{mN1}
\lim_{t\ri+\iy}m_t[\cN]&=&1\eqn
\end{theo}
Thus to find a element of $\cM$ with an important probability, one should 
pick up the value of $X_t$ for sufficiently large times $t$.\par

A crucial ingredient of the proof of this convergence are the Gibbs measures associated to the potentials $U_\delta$.
They are defined as the probability measures $\mu_{\beta,\delta}$ given for
 any $\beta\geq0$ by
\bqn{mub}
\mu_{\beta,\delta}(dx)&\df&\frac{\exp(-\beta U_\delta(x))}{Z_{\beta,\delta}}\, \lambda(dx)\eqn
where $Z_{\beta,\delta}\df \int \exp(-\beta U_\delta(x))\, \lambda(dx)$ is the normalizing factor.
\\
Indeed we will show that $\cL(X_t)$ and $\mu_{\beta_t,\delta_t}$
become closer and closer as $t\geq0$ goes to infinity  in the sense
of  $\LL^2$ variation. 

By uniform continuity of $\kappa$ we can easily prove that $U_\delta$ converges uniformly to $U$ as $\delta\to 0$. As a consequence, for any neighbourhood $\cN$ of $\cM$, $\mu_{\beta,\delta}(\cN)$ converges to $1$ as $\beta\to\infty$, uniformly in $\delta\le \delta_0$ for some $\delta_0>0$ depending on $\cN$. All this will prove the theorem.

The main difference in the method between the present work and~\cite{Arnaudon} concerns the definition of the jumps. Instead of following the geodesic from the current position to a realization of $\nu$, the process jumps to $\phi_{\delta_t}(\beta_t\alpha_t, X_{t-}, Y_{N_t^{(\alpha)}})$. The calculations are much easier, and this allows to consider more general cost functions $\kappa$. The drawback is that the implementation may be more complicated.
 
 The cost functions $\kappa$ are only assumed to be continuous. So this includes distances at the power $p$, $\rho^p$, which lead to $p$-means, for all $p\in (0,\infty)$. Notice the case $p\in(0,1)$ has never been considered in previous works. This also includes lengths for Finsler metrics, all divergences  for parametric statistical models (Kullback-Leibler, Jeffrey, Chernoff, Hellinger...).
\sm
\par\sm
The paper is constructed on the following plan.
In  Section 2 we obtain an estimate of $L_{\alpha,\beta,\delta}^\ast\un$ where $L_{\alpha,\beta,\delta}^\ast$ is the adjoint of $L_{\alpha,\beta,\delta}$ in $L^2(\mu_{\beta,\delta})$, $L_{\alpha,\beta,\delta}$ is the generator of the process $X_t$ described above but with constant $\alpha, \beta, \delta$. The proof of Theorem~\ref{t1} is given in Section 3.
 For this proof, the estimate of Section 2 is crucial to see how close is the instantaneous invariant measure associated to
  the algorithm at large times $t\geq0$ to the Gibbs measures associated to
  the potential $U_{\delta_t}$ and to the inverse temperature $\beta_t^{-1}$.

\section{Regularity issues}

\par\me
Let us consider a general probability measure $\nu$ on $M$.
For any $\alpha>0$, $\beta\geq 0$ and $\delta>0$, we are interested into the generator $L_{\alpha,\beta,\delta }$
defined for $f$ from $\cC^2(M)$ via
\bqn{Lab}
\fo x\in M,\qquad L_{\alpha,\beta,\delta}[f](x)&\df & \frac12\tr f(x)+\frac1{\alpha}\int \left(f(\phi_\delta(\beta\alpha,x,y))-f(x)\right)\,\nu(dy)\eqn
\par\sm
 We will prove in Section 3 that   $\cL(X_t)$ gets  closer and closer to   the Gibbs distribution $\mu_{\beta_t,\delta_t}$ as $t\to\infty$. Since for large $\beta\geq 0$, $\mu_{\beta,\delta}$ concentrates
 around $\cM$ uniformly in $\delta$ sufficiently small, this will be sufficient to  establish Theorem~\ref{t1}. 
 \par
 The remaining part of this section is devoted to a quantification of what separates $\mu_{\beta,\delta}$ from
 being an invariant probability of $L_{\alpha,\beta,\delta}$, for
 $\alpha>0$, $\delta>0$ and $\beta\geq 0$. As it will become clearer in the next section, a practical
 way to measure this discrepancy is through the evaluation 
 of $\mu_{\beta,\delta}[(L^*_{\alpha,\beta,\delta}[\un])^2]$, where $L^*_{\alpha,\beta,\delta}$ is the dual operator
 of $L_{\alpha,\beta,\delta}$ in $\LL^2(\mu_{\beta,\delta})$ and where $\un$ is the constant function
 taking the value 1. Indeed, it can be seen that $L^*_{\alpha,\beta,\delta}[\un]=0$ in  $\LL^2(\mu_{\beta,\delta})$
 if and only if $\mu_{\beta,\delta}$ is invariant for $L_{\alpha,\beta,\delta}$. Before being more precise about the definition of $L^*_{\alpha,\beta,\delta}$, we need an elementary result, where we will use the following notations:  for any $s\in\RR$, $T_{\delta,y,s}$ is the operator acting on measurable functions $f$ defined on $M$ via
 \bqn{Tys}\fo x\in M,\qquad T_{\delta,y,s}f(x)&\df& f(\phi_\delta(s,x,y))\eqn
 \begin{lem}\label{gf}
 For any $y\in M$, any $s\in [0,1)$ and any measurable and bounded functions $f,g$, we have
 \bq
 \int_M g(x)T_{\delta,y,s}f(x)\,\lambda(dx)&=&\int_M f(z)T_{\delta,y,-s}g(z)|J\phi_\delta(-s, \cdot, y)|(z)\,\lambda(dz)\eq
 where $dx$ and $dz$ denote Lebesgue measure on $M$ and $J\phi_\delta(-s, \cdot, y)(z)$ is the determinant at $z$ of the Jacobian matrix of $\phi_\delta(-s, \cdot, y)$.  
 \end{lem}
 \proof 
 Just make the change of variable $z=\phi_\delta(s,x,y)$ in the first integral, which yields $x=\phi_\delta(-s,z,y)$.\wwtbp
\par
This lemma has for consequence the next result,  where $\cD$ is the subspace of $\LL^2(\lambda)$
consisting of functions whose second derivatives in the distribution sense belongs
to $\LL^2(\lambda)$ (or equivalently to $\LL^2(\mu_{\beta,\delta})$ for any $\beta\geq 0$ and $\delta>0$).
 \begin{lem}
 For $\alpha>0$, $\beta\geq 0$ and $\delta>0$,
  the domain of the
 maximal extension of $L_{\alpha,\beta,\delta}$ on $\LL^2(\mu_{\beta,\delta})$ is $\cD$.
 Furthermore the domain of its dual operator $L_{\alpha,\beta,\delta}^*$ in  $\LL^2(\mu_{\beta,\delta})$
 is also $\cD$ and we have for any $f\in \cD$, 
 \bq
\lefteqn{ L_{\alpha,\beta,\delta}^*f}\\
 &=&
\frac12 \exp(\beta U_\delta)\Delta[\exp(-\beta U_\delta)f]+\frac{\exp(\beta U_\delta)}{\alpha}\int T_{\delta,y,-\alpha\beta}[\exp(-\beta U_\delta)f ]|J\phi_\delta(-\alpha\beta, \cdot, y)|\, \nu(dy)-\frac{f}{\alpha}\eq
 \end{lem}
 \proof
With the previous definitions, we can write for any $\alpha>0$, $\beta\geq 0$ and $\delta>0$,
 \bq
 L_{\alpha,\beta,\delta}&=&\frac12\Delta+\frac1{\alpha}\int T_{\delta,y,\alpha\beta}\,\nu(dy)-\frac{I}{\alpha}\eq
 where $I$ is the identity operator.
Note furthermore that the identity operator is bounded from $\LL^2(\lambda)$
to $\LL^2(\mu_{\beta,\delta})$ and conversely.
Thus to get the first assertion, it is sufficient to show that
$\int T_{\delta,y,\alpha\beta}\,\nu(dy)$ is bounded from $\LL^2(\lambda)$ to itself, 
or even only that $\lVe T_{\delta,y,\alpha\beta}\rVe_{\LL^2(\lambda)\righttoleftarrow}$
is uniformly bounded in $y\in M$.
To see that this is true, consider a bounded and measurable function $f$ and assume that
$\alpha\beta\ge 0$.
Since $(T_{\alpha\beta}f)^2=T_{\alpha\beta}f^2$, we can apply Lemma \ref{gf}
with $s=\alpha\beta$, $g=\un$ and $f$ replaced by $f^2$ to get
that
\bq
\int (T_{\delta,y,\alpha\beta}f)^2(x)\,\lambda(dx)&=&
\int f^2(z)T_{\delta,y,-\alpha\beta}\un |J\phi_\delta(-\alpha\beta, \cdot, y)|(z)\,\lambda(dz)\\
&\leq &J_{\delta,\infty}\int f^2\,d\lambda
\eq
with $\di J_{\delta,\infty}=\sup_{z,y\,\in M}|J\phi_\delta(-\alpha\beta, \cdot, y)|(z)$.
This quantity is  finite, since $\kappa_{\delta}(\cdot,\cdot)$ belongs to the class $\cC^{\iy,0}$,
due to its definition by convolution with a smooth kernel.
 Next to see that for any $f,g\in\cC^2(M)$, 
 \bqn{Ladjoint}
 \int g L_{\alpha,\beta,\delta}f\,d\mu_{\beta,\delta}&=&  \int f L_{\alpha,\beta,\delta}^*g\,d\mu_{\beta,\delta}\eqn
 where $L_{\alpha,\beta,\delta}^*$ is the operator defined in the statement of the lemma,
 we note that, on one hand,
 \bq
  \int g \Delta f\,d\mu_{\beta,\delta}&=&Z_{\beta,\delta}^{-1}\int \exp(-\beta U_\delta)g\Delta f\,d\lambda\\
  &=& \int f \exp(\beta U_\delta)\Delta [\exp(-\beta U_\delta)g]\,d\mu_{\beta,\delta}\eq
  and on the other hand, for any $y\in M$,
  \bq
  \int gT_{\delta,y,\alpha\beta}f \,d\mu_{\beta,\delta}&=&
  Z_{\beta,\delta}^{-1}\int \exp(-\beta U_\delta)gT_{\delta,y,\alpha\beta}f\,d\lambda\\
  &=& Z_{\beta,\delta}^{-1}\int_M f T_{\delta,y,-\alpha\beta}\left(\exp(-\beta U_\delta)g\right)|J\phi_\delta(-\alpha\beta,\cdot,y)|(x)\,\lambda(dx)
  \eq
  by Lemma \ref{gf}.
  After an additional  integration with respect to $\nu(dy)$,
  (\ref{Ladjoint})  follows without difficulty.
 To conclude, it is sufficient to see that for any $f\in\LL^2(\mu_{\beta,\delta})$,
 $L^*_{\alpha,\beta,\delta}f\in \LL^2(\mu_{\beta,\delta})$ (where $L^*_{\alpha,\beta,\delta}f$ is first interpreted as a distribution) if and only if $f\in\cD$.
 This is done by adapting the arguments given in the first part of the proof, in particular we get
 that
 \bq
 \lVe 
\frac{\exp(\beta U_\delta)}{\alpha}\int    T_{\delta,y,-\alpha\beta}[\exp(-\beta U_\delta)\,\cdot\, ]|J\phi_\delta(-\alpha\beta,\cdot,y)|\, \nu(dy)\rVe^2_{\LL^2(\lambda)\righttoleftarrow}&\leq &
\frac{J_{\delta,\infty}^2\exp(2\beta\osc(U_\delta))}{\alpha^2}.
 \eq
 \wwtbp
 
For any $\alpha>0$ and $\beta\geq 0$, 
denote $\eta=\alpha\beta$.
As a consequence of the previous lemma, we get that for any $x\in M$,
\bqn{Lstarun}
\nonumber L_{\alpha,\beta,\delta}^*\un(x)&=&
\frac12 \exp(\beta U_\delta(x))\Delta\exp(-\beta U_\delta(x))-\frac{1}{\alpha}\\&&+
\nonumber \frac{\exp(\beta U_\delta(x))}{\alpha}\int  T_{\delta,y,-\eta}[\exp(-\beta U_\delta)](x)|J\phi_\delta(-\eta,\cdot,y)|(x)\, \nu(dy)\\
\nonumber &=&\frac{\beta^2}2(|\na U_\delta|(x))^2-\frac{\beta}2\Delta U_\delta(x)-\frac{1}{\alpha}\\
&&+\frac{1}{\alpha}\int_M 
\exp(\beta  [U_\delta(x)-U_\delta(\phi_\delta(-\eta,x,y))])|J\phi_\delta(-\eta,\cdot,y)|(x)
  \, \nu(dy)\eqn
  It appears that $L_{\alpha,\beta,\delta}^*\un$ is continuous. The next result
  evaluates the uniform norm of this function.
  
  \begin{pro}\label{mbl0}
  There exists  a constant $C>0$, depending on $M$ and $\lVe \kappa\rVe_\iy$, such that
  for any $\beta\geq 1$,  $\delta\in (0,1]$ and $\alpha\in (0,\delta^2/(2\beta^2))$ we have
  \bq
  \lVe L_{\alpha,\beta,\delta}^*\un\rVe_\iy &\leq &C\alpha\beta^4\delta^{-4}\eq
  \end{pro}
  \proof
  In view of the expression of $ L_{\alpha,\beta,\delta}^*\un(x)$ given before
  the statement of the proposition, we want to estimate for any fixed $x\in\TT$,
  the quantity
  \bq
  \int_M
\exp(\beta  [U_\delta(x)-U_\delta(\phi_\delta(-\eta,x,y))])|J\phi_\delta(-\eta,\cdot,y)|(x)
  \, \nu(dy).\eq
  Consider the function
  $$
  \psi(s)=U_\delta(x)-U_\delta(\phi_\delta(s,x,y)).
  $$
  It has derivative
  $$
  \psi'(s)=\frac12 \langle \na U_\delta, \na \kappa_\delta (\cdot, y)\rangle(\phi_\delta(s,x,y))
  $$
  and second derivative 
  $$
  \psi''(s)=-\frac14\Hess\  U_\delta \left(\na \kappa_\delta (\cdot, y)(\phi_\delta(s,x,y)),\na \kappa_\delta (\cdot, y)(\phi_\delta(s,x,y))\right)-\frac14\left\langle \na U_\delta, \na_{\na \kappa_\delta (\cdot, y)(\phi_\delta(s,x,y))}\na \kappa_\delta(\cdot, y)\right\rangle.
  $$
 For any $\eta=\alpha\beta$, there exists $s\in[0,\eta]$ such that 
 $$
 \psi(-\eta)=\psi(0)- \eta\psi'(0)+\frac{\eta^2}{2}\psi''(-s).
 $$
 This yields 
 \begin{align*}
 &\beta\left(U_\delta(x)-U_\delta(\phi(-\eta,x,y))\right)=\frac{-\beta\eta}2 \langle \na U_\delta,\na \kappa_\delta(\cdot,y)\rangle(x)\\&-\frac{\beta\eta^2}8\left(\Hess\  U_\delta \left(\na \kappa_\delta (\cdot, y)(\phi_\delta(-s,x,y)),\na \kappa (\cdot, y)(\phi_\delta(-s,x,y))\right)+\left\langle \na U_\delta, \na_{\na \kappa (\cdot, y)(\phi_\delta(-s,x,y))}\na \kappa(\cdot, y)\right\rangle\right)
 \end{align*}

  Observe that for any $a,b\in\RR$, we can find $u,v\in(0,1)$ such that
  \begin{equation}
  \label{expaplusb}
  \exp(a+b)=(1+a+a^2\exp(ua)/2)(1+b\exp(vb)).
  \end{equation} Apply this equality with 
  $$
  a=\frac{-\beta\eta}2 \langle \na U_\delta,\na \kappa_\delta(\cdot,y)\rangle(x)
  $$ 
  and 
  $$
  b=-\frac{\beta\eta^2}8\left(\Hess\  U_\delta \left(\na \kappa_\delta (\cdot, y)(\phi_\delta(s,x,y)),\na \kappa_\delta (\cdot, y)(\phi_\delta(-s,x,y))\right)+\left\langle \na U_\delta, \na_{\na \kappa_\delta (\cdot, y)(\phi_\delta(-s,x,y))}\na \kappa_\delta(\cdot, y)\right\rangle\right).
  $$
  Using the bounds 
  \begin{equation}
  \label{gradlnp}
  \fo \delta>0,\qquad
  \|\nabla \ln p(\delta,\cdot,y)(x)\|\le\frac{C'}{\delta}\quad\hbox{and}\quad \|\nabla d \ln p(\delta,\cdot,y)(x)\|\le\frac{C'}{\delta^2}
  \end{equation}
  for some $C'>0$ (see e.g. \cite{Hsu:99}), writing
  $$
  \nabla \kappa_\delta(\cdot,y)(x)=\int_M\nabla \ln p(\delta,\cdot,z)p(\delta,x,z)\kappa(z,y)\,\lambda(dz),
  $$
  we get
  \begin{equation}
  \label{gradlnkappa} \fo \delta>0,\qquad
  \|\nabla \kappa_\delta(\cdot,y)(x)\|\le\frac{C}{\delta}\quad\hbox{and}\quad \|\nabla d\kappa_\delta(\cdot,y)(x)\|\le\frac{C}{\delta^2}
  \end{equation}
  with $C=2C'\|\kappa\|_\infty$,
  together with
  \begin{equation}
  \label{gradlnU}
  \|\nabla  U_\delta(x)\|\le\frac{C}{\delta}\quad\hbox{and}\quad \|\nabla d U_\delta(x)\|\le\frac{C}{\delta^2}.
  \end{equation}
  It follows that $\lve a\rve =\cO(\alpha\beta^2\delta^{-2})$ and $\lve b\rve =\cO(\alpha^2\beta^3\delta^{-4})$, so
  in conjunction with the assumption $\alpha\beta^2\delta^{-2}\leq 1/2$, we can write with~\eqref{expaplusb} that
  \bqn{expUU}
   \exp(\beta [U_\delta(x)-U_\delta(\phi_\delta(-\eta,x,y))])&=&1- \frac{\beta \eta}{2} \langle \na U_\delta,\na \kappa_\delta(\cdot,y)\rangle(x)+\cO(\alpha^2\beta^4\delta^{-4}).\eqn
   Integrating this expression, we get that
   \bq
  \lefteqn{  \int_M 
\exp(\beta  [U_\delta(x)-U_\delta(\phi(-\eta,x,y))])|J\phi_\delta(-\eta,\cdot,y)|(x)
  \, \nu(dy)}\\&=&
  \int_M 
|J\phi_\delta(-\eta,\cdot,y)|(x)
  \, \nu(dy)\\&-&
  \frac{\beta \eta}2  \int_M\langle \na U_\delta,\na \kappa_\delta(\cdot,y)\rangle(x)|J\phi_\delta(-\eta,\cdot,y)|(x)\, \nu(dy)
  +\cO(\alpha^2\beta^4\delta^{-4})\\
\eq
where we used the fact that $|J\phi_\delta(-\eta,\cdot,y)|(x)$ is uniformly bounded (see~\eqref{devJac} below).
 We can now return to (\ref{Lstarun}) and  we obtain that for any $x\in M$,
 \bq
 L_{\alpha,\beta,\delta}^*\un(x)&=&
\frac{\beta^2}2|\na U_\delta(x)|^2-\frac{\beta}{2}\Delta U_\delta(x)+\frac1\alpha \int_M\left(|J\phi_\delta(-\eta,\cdot,y)|(x)-1\right)\,\nu(dy)\\
&&-\frac{\beta^2}2  \int_M\langle \na U_\delta,\na \kappa(\cdot,y)\rangle(x)|J\phi_\delta(-\eta,\cdot,y)|(x)\, \nu(dy)+\cO(\alpha\beta^4\delta^{-4}).
\eq
Note that for $\eta/\delta\geq 0$ small enough (up to a universal factor, less than the injection radius of $M$), we have
$$
\phi_\delta(-\eta,x,y)=\exp_x\left(\frac{\eta}2\na\kappa_\delta(\cdot,y)(x)+\eta^2 \int_0^1\left(\log_x\circ\phi_\delta\right)''(-s\eta,x,y)(1-s)\,ds)\right)
$$

where $\log_x$ is the inverse function of $\exp_x$ and $\left(\log_x\circ\phi_\delta\right)''$ is the second derivative in the first variable.
From this equality,  in conjunction with \eqref{phiprime} and~\eqref{gradlnkappa},
we get 
\begin{equation}
\label{devJac}
|J\phi_\delta(-\eta,\cdot,y)|(x)=1+\frac{\eta}2\Delta \kappa_\delta(\cdot,y)(x)+\cO(\eta^2\delta^{-4})
\end{equation}
first for $\alpha\beta/\delta^2$ small enough and next by a compactness argument for all $\alpha,\beta,\delta$ in the range described in the statement of Proposition
\ref{mbl0}.
It  also appears that $|J\phi_\delta(-\eta,\cdot,y)|(x)$ is uniformly bounded when $\alpha\beta^2\delta^{-2}\le\frac12$.
This yields 
\bq
\frac1\alpha \int_M\left(|J\phi_\delta(-\eta,\cdot,y)|(x)-1\right)\,\nu(dy)
&=&\frac{\beta}2\int_M\Delta \kappa_\delta(\cdot,y)(x)\,\nu(dy)+\cO(\alpha\beta^2\delta^{-4})\\
&=&\frac{\beta}{2}\Delta U_\delta(x)+\cO(\alpha\beta^2\delta^{-4}).
\eq
Notice the first term in the right cancels with the second in the right of~\eqref{Lstarun}.
We also have
\begin{align*}
\lefteqn{-\frac{\beta^2}2  \int_M\langle \na U_\delta,\na \kappa_\delta(\cdot,y)\rangle(x)|J\phi_\delta(-\eta,\cdot,y)|(x)\, \nu(dy)}\\&=-\frac{\beta^2}2\left\langle\na U_\delta,\int_M\na \kappa_\delta(\cdot,y)(x)\, \nu(dy)\right\rangle+\cO(\alpha\beta^3\delta^{-4})\\
&=\frac{-\beta^2}2|\na U_\delta(x)|^2+\cO(\alpha\beta^3\delta^{-4}).
\end{align*} 
Here the first term in the right cancels with the first term in the right of~\eqref{Lstarun}.
The bound announced in the lemma follows at once.
\wwtbp
In particular, under the hypotheses of the previous
proposition we get
\bqn{mbl}
\sqrt{\mu_{\beta,\delta}[(L^*_{\alpha,\beta,\delta}\un)^2]}&\leq &C\alpha\beta^4\delta^{-4}\eqn
The l.h.s.\ will be used in the next section,  when $\alpha\beta^4\delta^{-4}$
is small, as a discrepancy for the fact that $\mu_{\beta,\delta}$
is not necessarily an invariant measure for $L_{\alpha,\beta,\delta}$.
\par\me

\section{Proof of Theorem 1}\label{proof}

This is the main part of the paper: we are going to prove Theorem \ref{t1}
by the investigation of the evolution of a $\LL^2$ type functional.\par\me
On $M$ consider the algorithm $X\df(X_t)_{t\geq 0}$ described in the introduction.
For the time being, 
 the schemes $\alpha\st \RR_+\ri \RR_+^*$,
$\beta\st \RR_+\ri \RR_+$ and $\delta\st \RR_+\ri \RR_+^\ast$ are assumed to be  continuously differentiable.
Only later on, in Proposition  \ref{condab}, will we present the conditions insuring the wanted convergence (\ref{mN1}).
On the initial distribution $m_0$, the last ingredient necessary to specify
the law of $X$, no hypothesis is made.
We also denote  $m_t$ the law of $X_t$, for any $t>0$. We have that $m_t$ admits a $\cC^1$ density with respect to $\lambda$,
which is equally written $m_t$ (for a proof we refer to the appendix of~\cite{Arnaudon}).
As it was mentioned in the previous section, we want to compare these temporal marginal laws
with the corresponding instantaneous Gibbs measures, which were defined in (\ref{mub})
with respect to the potentials $U_\delta$ given in (\ref{Us}).
A convenient way to quantify this discrepancy is to consider the variance of the density of 
$m_t$ with respect to $\mu_{\beta_t,\delta_t}$ under the probability measure $\mu_{\beta_t,\delta_t}$:
\bqn{It}
\fo t>0,\qquad
I_t&\df& \int \lt(\frac{m_t}{\mu_{\beta_t,\delta_t}}-1\rt)^2\,d\mu_{\beta_t,\delta_t}\eqn
Our goal here is to derive a differential inequality satisfied by this quantity,
which implies its convergence to zero under appropriate conditions on the
schemes $\alpha$ and $\beta$. 
More precisely, our purpose is to obtain:
\begin{pro}\label{Iprime}
There exists two constants $c_1,\,c_2>0$  such that
for any $t>0$ with $\beta_t\geq 1$ and $\alpha_t\beta_t^2\delta_t^{-2}\leq 1/2$, we have
\bq
I_t'&\leq & -c_1\left[(\beta_t\delta_t^{-1})^{2-5m}\exp(-b(U)\beta_t)-\alpha_t\beta_t^3\delta_t^{-3}-\lve \beta_t'\rve-\beta_t\delta_t^{-2}\lve \delta_t'\rve\right]I_t\\&+&c_2\left[\alpha_t\beta_t^4\delta_t^{-4}+\lve \beta_t'\rve +\beta_t\delta_t^{-2}\lve \delta_t'\rve\right]\sqrt{I_t}\eq
where $b(U)$ was defined in (\ref{bU}).
\end{pro}
\proof
At least formally, there is no difficulty to differentiate the quantity $I_t$
with respect to the time $t>0$. 
For a  rigorous justification of the following computations, we refer to the appendix of 
\cite{Arnaudon}, where the regularity of the temporal marginal laws in presence of jumps is discussed in detail
(it is written in the situation considered there of the circle but can be extended to compact manifolds).
Thus we get at any time $t>0$,
\bq
I'_t&=&2\int \lt(\frac{m_t}{\mu_{\beta_t,\delta_t}}-1\rt)\frac{\pa_t m_t}{\mu_{\beta_t,\delta_t}}\,d\mu_{\beta_t,\delta_t}
-2\int \lt(\frac{m_t}{\mu_{\beta_t,\delta_t}}-1\rt)\frac{m_t}{\mu_{\beta_t,\delta_t}}\pa_t \ln(\mu_{\beta_t,\delta_t})\,d\mu_{\beta_t,\delta_t}\\&&+\int \lt(\frac{m_t}{\mu_{\beta_t,\delta_t}}-1\rt)^2\pa_t \ln(\mu_{\beta_t,\delta_t})\,d\mu_{\beta_t,\delta_t}\\
&=&
2\int \lt(\frac{m_t}{\mu_{\beta_t,\delta_t}}-1\rt)\pa_t m_t\,d\lambda
-\int \lt(\frac{m_t}{\mu_{\beta_t,\delta_t}}-1\rt)^2\pa_t \ln(\mu_{\beta_t,\delta_t})\,d\mu_{\beta_t,\delta_t}
\\&-&2\int \lt(\frac{m_t}{\mu_{\beta_t,\delta_t}}-1\rt)\pa_t \ln(\mu_{\beta_t,\delta_t})\,d\mu_{\beta_t,\delta_t}\\
&\leq & 2\int \lt(\frac{m_t}{\mu_{\beta_t,\delta_t}}-1\rt)\pa_t m_t\,d\lambda
+\lVe \pa_t \ln(\mu_{\beta_t,\delta_t})\rVe_{\iy}\lt(
\int \lt(\frac{m_t}{\mu_{\beta_t,\delta_t}}-1\rt)^2\,d\mu_{\beta_t,\delta_t}+
2\int \lve\frac{m_t}{\mu_{\beta_t,\delta_t}}-1\rve\,d\mu_{\beta_t,\delta_t}\rt)\\
&\leq &  2\int \lt(\frac{m_t}{\mu_{\beta_t,\delta_t}}-1\rt)\pa_t m_t\,d\lambda
+\lVe \pa_t \ln(\mu_{\beta_t,\delta_t})\rVe_{\iy}\lt(
I_t+2\sqrt{I_t}\rt)
\eq
where we used the Cauchy-Schwarz inequality.
The last term is easy to deal with:
\begin{lem}\label{palnmu}
There exists $C_0\geq 0$, depending on $\kappa$, such that for any $t\geq 0$, we have
\bq
\lVe \pa_t \ln(\mu_{\beta_t,\delta_t})\rVe_{\iy}&\leq & C_0\left(\lve \beta_t'\rve +\beta_t\lve \delta_t'\rve \delta_t^{-2}\right).\eq
\end{lem}
\proof
Since for any $t\geq 0$ we have
\bq
\fo x\in M,\qquad\ln(\mu_{\beta_t,\delta_t})(x)&=&- \beta_tU_{\delta_t} (x)-\ln\lt(\int \exp(-\beta_t U_{\delta_t}(y))\,\lambda(dy)\rt)\eq
it appears that $\fo x\in M,$
\bq
\lefteqn{\pa_t\ln(\mu_{\beta_t})(x)}\\&=&\beta_t'\int U_{\delta_t}(y)-U_{\delta_t}(x)\, \mu_{\beta_t,\delta_t}(dy)
\\&&+\beta_t\delta_t'\int\!\int\!\int\left(p(\delta_t,y,z)\partial_\delta \ln p(\delta_t,y,z)-p(\delta_t,x,z)\partial_\delta \ln p(\delta_t,x,z)\right)\kappa(z,v)\nu(dv)\mu_{\beta_t,\delta_t}(dy)\lambda(dz)
\eq
so that
\bq
\lVe \pa_t \ln(\mu_{\beta_t})\rVe_{\iy}&\leq & \osc(U_{\delta_t})\lve \beta_t'\rve+2\beta_t|\delta_t'|\|\partial_\delta \ln p\|_\infty\cdot\|\kappa\|_\infty.\eq
Clearly  $\osc(U_{\delta_t})\leq 2\|\kappa\|_\infty$. To finish the proof we are left to use the bound
\begin{equation}
\label{partialnp}
|\partial_\delta \ln p(\delta,x,y)|\le \frac{C''}{\delta^2}
\end{equation}
for some $C''>0$ (see e.g.~\cite{Hsu:99}).
\wwtbp\par
Denote for any $t>0$, $f_t\df m_t/\mu_{\beta_t,\delta_t}$.
If this function was to be $\cC^2$, we would get, by the martingale problem satisfied by the
law of $X$, that
\bq
\int \lt(\frac{m_t}{\mu_{\beta_t,\delta_t}}-1\rt)\pa_t m_t\,d\lambda&=&
\int L_{\alpha_t,\beta_t,\delta_t}\lt[ 
f_t-1\rt]
\, dm_t\\
&=&\int L_{\alpha_t,\beta_t,\delta_t}\lt[ 
f_t-1\rt] f_t
\, d\mu_{\beta_t,\delta_t}\eq
where $ L_{\alpha_t,\beta_t,\delta_t}$, described in the previous section, is the instantaneous generator at time $t\geq 0$ of
$X$. The interest of the estimate (\ref{mbl}) comes from the decomposition of the previous term
into
\bq
\lefteqn{\int L_{\alpha_t,\beta_t,\delta_t}\lt[ 
f_t-1\rt] (f_t-1)
\, d\mu_{\beta_t,\delta_t}+\int L_{\alpha_t,\beta_t,\delta_t}\lt[ 
f_t-1\rt] \, d\mu_{\beta_t,\delta_t}}\\
&=&
\int L_{\alpha_t,\beta_t,\delta_t}\lt[ 
f_t-1\rt] (f_t-1)
\, d\mu_{\beta_t,\delta_t}+\int(
f_t-1)L_{\alpha_t,\beta_t,\delta_t}^*[\un] \, d\mu_{\beta_t,\delta_t}\\
&\leq &
\int L_{\alpha_t,\beta_t,\delta_t}\lt[ 
f_t-1\rt] (f_t-1)
\, d\mu_{\beta_t,\delta_t}+\sqrt{I_t}\sqrt{\mu_{\beta_t,\delta_t}[(L_{\alpha_t,\beta_t,\delta_t}^*[\un])^2]}\eq
It follows from these bounds that to prove Proposition \ref{Iprime}, it 
remains to treat the first term in the above r.h.s. A first step is:
\begin{lem}\label{Lffm}
There exists a constant $c_3>0$, such that
for any $\alpha>0$ and $\beta\geq 1$ such that $\alpha\beta^2\delta^{-2}\leq 1/2$, we have,
for any $f\in\cC^2(M)$,
\bq\int L_{\alpha,\beta,\delta}\lt[ 
f-1\rt] (f-1)
\, d\mu_{\beta,\delta}&\leq &-\lt(\frac{
1}2-c_3\alpha\beta^3\delta^{-3}\rt)\int (|\na f|)^2\, d\mu_{\beta,\delta}
+c_3\alpha\beta^3\delta^{-3}\int (f-1)^2\, d\mu_{\beta,\delta}
\eq
\end{lem}
\proof
For any $\alpha>0$, $\beta\geq 0$ and $\delta>0$, we begin by decomposing 
the generator $L_{\alpha,\beta,\delta}$ into
\bqn{LLR}
L_{\alpha,\beta,\delta}&=&L_{\beta,\delta}+R_{\alpha,\beta,\delta}\eqn
where 
\bqn{Lb}
\qquad L_{\beta,\delta}\,\cdot&\df& \frac{1}{2}(\tr\cdot -\beta \lan \na U_\delta,\na\cdot \ran)\eqn
and where $R_{\alpha,\beta,\delta}$ is the remaining operator.
An immediate integration by parts leads to 
\begin{equation}
\label{iip}
\begin{split}
\int L_{\beta,\delta}\lt[ 
f-1\rt] (f-1)
\, d\mu_{\beta,\delta}= &-\frac{
1}{2}\int |\na (f-1)|^2\, d\mu_{\beta,\delta}\\
=&-\frac{
1}{2}\int |\na f|^2\, d\mu_{\beta,\delta}
\end{split}
\end{equation}
Thus our main task is to find a 
constant $c_3>0$, such that
for any $\alpha>0$, $\beta\geq 1$ and $\delta>0$ with $\alpha\beta^2\delta^{-2}\leq 1/2$, we have,
for any $f\in\cC^2(M)$,
\bqn{borned}
\lve\int R_{\alpha,\beta,\delta}\lt[ 
f-1\rt] (f-1)
\, d\mu_{\beta,\delta}\rve&\leq &
{c_3\alpha\beta^3\delta^{-3}}\lt(\int |\na f|^2\, d\mu_{\beta,\delta}+\int(f-1)^2\, d\mu_{\beta,\delta}\rt)
\eqn
By definition, we have for any $f\in\cC^2(M)$ (but what follows is valid for $f\in\cC^1(M)$),
\bq
\fo x\in M,\qquad
R_{\alpha,\beta,\delta}[f](x)&=& \frac1\alpha\int f(\phi_\delta(\alpha\beta, x,y))-f(x)\,\nu(dy)+\frac\beta2\langle \na U_\delta(x),\na f(x)\rangle\eq
To evaluate this quantity, on one hand, recall that we have for any $x\in M$,
\bq
\na U_\delta(x)&=&\int_{M} \na \kappa_\delta(\cdot, y)(x)\,\nu(dy)\eq
and on the other hand, write that for any $x,y\in M$,
\bq
f(\phi_\delta(\alpha\beta,x,y))-f(x)&=&\alpha\beta\int_0^1\left\langle \na f(\phi_\delta(\alpha\beta u,x,y )),-\frac12\na\kappa_\delta(\cdot,y)(\phi_\delta(\alpha\beta u,x,y ))\right\rangle\,du\eq
It follows that 
\bq
\lefteqn{\int R_{\alpha,\beta,\delta}\lt[ 
f-1\rt] (f-1)
\, d\mu_{\beta,\delta}}\\&=&
\frac{\beta}2\int_0^1du\int \nu(dy)\int\mu_{\beta,\delta}(dx) \left(\langle \na f, \na \kappa_\delta(\cdot, y)(x)\rangle-\langle \na f,\na\kappa_\delta(\cdot, y)(\phi_\delta(\alpha\beta u,x,y))\right)(f(x)-1)\\
&=&\frac{\beta}2\int_0^1du\int \nu(dy)\int \lambda(dx)\langle \na f, \na \kappa_\delta(\cdot, y)(x)\rangle\Biggl[(f(x)-1)\mu_{\beta,\delta}(x)\\
&-&\left(f(\phi_\delta(-\alpha\beta u,x,y))-1\right)\mu_{\beta,\delta}(\phi_\delta(-\alpha\beta u,x,y))|J\phi_\delta(-\alpha\beta u,\cdot,y)|(x)\Biggr]\eq
where we used the change of variable $z\mapsto \phi_\delta(-\alpha\beta u,z,y)$ for the second term in the right. So
\bq
\lefteqn{\int R_{\alpha,\beta,\delta}\lt[ 
f-1\rt] (f-1)
\, d\mu_{\beta,\delta}}\\
&=&
\frac{\beta}2\int_0^1du\int \nu(dy)\int\lambda(dx)\langle \na f, \na \kappa_\delta(\cdot, y)(x)\rangle I_\delta(\alpha\beta u,x,y)
\eq
where 
\begin{align*}
I_\delta(s,x,y)&=|J\phi_\delta(-s,\cdot,y)|(x)\left\{(f(x)-1)\mu_{\beta,\delta}(x)-(f(\phi_\delta(-s,x,y))-1)\mu_{\beta,\delta}(\phi_\delta(-s,x,y))\right\}\\
&+(f(x)-1)\mu_{\beta,\delta}(x)(1-|J\phi_\delta(-s,\cdot,y)|(x)).
\end{align*}
Write 
$$
\int R_{\alpha,\beta,\delta}\lt[ 
f-1\rt] (f-1)
\, d\mu_{\beta,\delta}=J_1+J_2
$$
where $J_1$ is the integral containing the term $(1-|J\phi_\delta(-s,\cdot,y)|(x))$.
From the validity of
$$
\phi_\delta(-s,x,y)=\exp_x\left(s\int_0^1(\log_x\circ\phi_\delta)'(-vs,x,y)\,dv\right)
$$
for $s$ small enough,
we get that for any $s$,
$$
\left|1-|J\phi_\delta(-s,\cdot,y)|(x)\right|\le c_4 s\delta^{-2}
$$
for some $c_4>0$, which yields
\begin{align*}
J_1&\le \frac12c_4\alpha\beta^2\delta^{-2}\int_0^1du\int \nu(dy)\int \mu_{\beta,\delta}(dx) \|\na_1\kappa_\delta\|_\infty |\na f|(x)\cdot |f(x)-1|\\
&\le \frac12c_4\alpha\beta^2\delta^{-2}\|\na_1\kappa_\delta\|_\infty \|\na f\|_{\LL^2(\mu_\beta)}\sqrt{\int(f-1)^2\, d\mu_\beta}\\
&\le \frac14c_4\alpha\beta^2\delta^{-3} \left(\|\na f\|_{\LL^2(\mu_{\beta,\delta})}^2+\int(f-1)^2\, d\mu_{\beta,\delta}\right)\\
&= c_5\alpha\beta^2\delta^{-3}\left(\|\na f\|_{\LL^2(\mu_{\beta,\delta})}^2+\int(f-1)^2\, d\mu_{\beta,\delta}\right)
\end{align*}
with  $c_5=\frac14c_4$, using again~\eqref{gradlnkappa}.
Moreover, we have 
\begin{align*}
J_2&=\frac{\beta}2\int_0^1du\int \nu(dy)\int \sqrt{\mu_{\beta,\delta}(x)}dx\langle \na f, \na \kappa_\delta(\cdot, y)\rangle(x)|J\phi_\delta(-\alpha\beta u,\cdot,y)|(x)(-\alpha\beta u)\int_0^1dv \\
&\Biggl(\langle \na f,-\frac12\na \kappa_\delta(\cdot,y)\rangle (\phi_\delta(-\alpha\beta uv,x,y))\sqrt{\mu_{\beta,\delta}(\phi_\delta(-\alpha\beta uv,x,y))}\sqrt{\frac{\mu_{\beta,\delta}(\phi_\delta(-\alpha\beta uv,x,y))}{\mu_{\beta,\delta}(x)}}\\&+f(\phi_\delta(-\alpha\beta uv,x,y))-1)\langle \na\ln\mu_{\beta,\delta},-\frac12\na\kappa_\delta(\cdot,y)\rangle(\phi_\delta(-\alpha\beta uv,x,y)) \\&\sqrt{\mu_{\beta,\delta}(\phi_\delta(-\alpha\beta uv,x,y))}\sqrt{\frac{\mu_{\beta,\delta}(\phi_\delta(-\alpha\beta uv,x,y))}{\mu_{\beta,\delta}(x)}} \Biggr).
\end{align*}
Recalling (\ref{devJac}), notice that $\di |J\phi_\delta(-\alpha\beta u,\cdot,y)|(x)$ and $\di \sqrt{\frac{\mu_{\beta,\delta}(\phi_\delta(-\alpha\beta uv,x,y))}{\mu_{\beta,\delta}(x)}}$ are uniformly bounded, since $\alpha\beta^2\delta^{-2}\leq 1/2$.
Putting the integral with respect to $v$ on the left, using Cauchy-Schwartz inequality for the integral in the right and making the change of variable $z=\phi_\delta(-\alpha\beta uv,x,y)$
we get 
\begin{equation}
\label{boundJ2}
\begin{split}
J_2&\le c_6\alpha\beta^2\|\na_1\kappa_\delta\|_\infty^2\|\na f\|_{\LL^2(\mu_{\beta,\delta})}\left(\|\na f\|_{\LL^2(\mu_{\beta,\delta})}+\beta\|\na_1\kappa_\delta\|_\infty\sqrt{\int(f-1)^2\,d\mu_{\beta,\delta}}\right)\\
&\le c_6\alpha\beta^2\|\na_1\kappa_\delta\|_\infty^2\left[\left(1+\frac12\beta\|\na_1\kappa_\delta\|_\infty\right)\|\na f\|_{\LL^2(\mu_{\beta,\delta})}^2+\frac12\beta\|\na_1\kappa_\delta\|_\infty^2\int(f-1)^2\,d\mu_{\beta,\delta}\right]
\end{split}
\end{equation}
where $c_6\geq 0$ is a constant  independent from $\alpha, \beta, \delta$. Up to a change of this constant, we obtain 
\begin{equation}
\label{bound2J2}
\begin{split}
J_2\le c_6\alpha\beta^3\delta^{-3}\left(\|\na f\|_{\LL^2(\mu_{\beta,\delta})}^2+\int(f-1)^2\, d\mu_{\beta,\delta}\right)
\end{split}
\end{equation}
where we used~\eqref{gradlnkappa}.

 So putting together~\eqref{iip} with the bounds for $J_1$ and $J_2$ 
 we get the wanted result.
\wwtbp
To conclude the proof of Proposition \ref{Iprime}, we must be able to compare,
for any $\beta\geq 0$ and any $f\in\cC^1(M)$, the energy $\mu_{\beta,\delta}[|\na f|^2]$
and the variance $\Var(f,\mu_{\beta,\delta})$. This task was already done by Holley, Kusuoka and Stroock \cite{MR995752}, let us recall their result:
\begin{pro}\label{HKS2}
Let $\wi U$ be a $\cC^1$ function on a compact Riemannian manifold $M$ of dimension $m\geq 1$.
Let $b(\wi U)\geq 0$ be the  associated constant  as in (\ref{bU}).
For any $\beta\geq 0$, consider the Gibbs measure  $\wi\mu_\beta$ 
given similarly to (\ref{mub}). Then there exists a constant $C_M>0$, depending  only on $M$,
such that the following Poincar\'e inequalities are satisfied:
\bq
\fo \beta\geq 0,\,\fo f\in\cC^1(M),\qquad \Var(f,\wi\mu_\beta)
&\leq & C_M[1\vee (\beta \lVe \nabla \wi U\rVe_\iy)]^{5m-2}\exp(b(\wi U)\beta)\wi\mu_\beta[\lve\na f\rve^2]\eq
\end{pro}
\par
We can now come back to the study of the evolution of the quantity $I_t=\Var(f_t,\mu_{\beta_t,\delta_t})$, for $t>0$. Indeed applying Lemma \ref{Lffm} and Proposition \ref{HKS2} with $\alpha=\alpha_t$,
$\beta=\beta_t$, $\delta=\delta_t$ and $f=f_t$, we get at any time $t>0$
such that $\beta_t\geq 1$, $\delta_t\in(0,1]$ and $\alpha_t\beta_t^2\delta_t^{-2}\leq 1/2$,
\bq
\lefteqn{\int L_{\alpha_t,\beta_t,\delta_t}\lt[ 
f_t-1\rt] (f_t-1)
\, d\mu_{\beta_t,\delta_t}}\\&\leq &-c_7(\beta_t\delta_t^{-1})^{2-5m}\exp(-b(U_{\delta_t})\beta_t)\lt(1
-2c_3\alpha_t\beta_t^3\delta_t^{-3}\rt)I_t+c_3\alpha_t\beta_t^3\delta_t^{-3}I_t\\
&\leq & 
-\left[c_7(\beta_t\delta_t^{-1})^{2-5m}\exp(-b(U_{\delta_t})\beta_t)-c_8\alpha_t\beta_t^3\delta_t^{-3}\right]I_t\eq
for some constants $c_7,c_8>0$.
\\
Taking into account Lemma \ref{palnmu}, the computations preceding Lemma  \ref{Lffm}
and (\ref{mbl}), one can find constants $c_1,c_2>0$ such that Proposition \ref{Iprime}
is satisfied. This achieves the proof of Proposition~\ref{Iprime}.
\wwtbp
\par\sm
This result leads immediately to conditions insuring the convergence toward 0
of the quantity $I_t$ for large times $t>0$:
\begin{pro}\label{condab}
Let 
$\alpha,\,\delta\st \RR_+\ri \RR_+^*$ and 
$\beta\st \RR_+\ri \RR_+$ be schemes as at the beginning of this section and assume:
\bq
\lim_{t\ri +\iy}\alpha_t&=&0\\
\lim_{t\ri +\iy}\beta_t&=&+\iy\\
\lim_{t\ri +\iy}\delta_t&=&0\\
\int_0^{+\iy}(1\vee(\beta_t\delta_t^{-1}))^{2-5m} \exp(-b(U_{\delta_t})\beta_t)\, dt &=&+\iy\eq
and that for large times $t>0$,
\bq
\max\{ \alpha_t\beta_t^4\delta_t^{-4},\, \lve \beta_t'\rve, \,\beta_t\delta_t^{-2}\lve \delta_t'\rve\}&\ll &  (\beta_t\delta_t^{-1})^{2-5m}\exp(-b(U_{\delta_t})\beta_t)\eq
Then we are assured of
\bq
\lim_{t\ri+\iy} I_t&=&0\eq
\end{pro}
\proof
The differential equation of Proposition \ref{Iprime} can be rewritten under the form
\bqn{Fee}
F_t'&\leq & -\eta_t F_t+\epsilon_t\eqn
where for any $t>0$,
\bq
F_t&\df& \sqrt{I_t}\\
\eta_t&\df & c_1((\beta_t\delta_t^{-1})^{2-5m}\exp(-b(U_{\delta_t})\beta_t)-\alpha_t\beta_t^3\delta_t^{-3}-\lve \beta_t'\rve-\beta_t\lve \delta_t'\rve\delta_t^{-2})/2\\
\epsilon_t&\df& c_2(\alpha_t\beta_t^4\delta_t^{-4}+\lve \beta_t'\rve+\beta_t\delta_t^{-2}\lve \delta_t'\rve)/2\eq
The assumptions of the above proposition 
imply that for $t\geq 0$ large enough,
$\beta_t\geq 1$ and $\alpha_t\beta_t^2\delta_t^{-2}\leq 1/2$.
This 
insures that there exists $T> 0$ such that (\ref{Fee})
is satisfied for any $t\geq T$ (and also $F_T<+\iy$). We deduce that for any $t\geq T$,
\bqn{Fprime}
F_t&\leq & F_T\exp\lt(-\int_T^t \eta_s\, ds\rt)
+\int_T^t \epsilon_s\exp\lt(-\int_s^t \eta_u\, du\rt)\, ds\eqn
It appears that $\lim_{t\ri +\iy} F_t =0$ as soon as
\bq
\int_T^{+\iy}\eta_s\, ds&=&+\iy\\
\lim_{t\ri +\iy} \epsilon_t/\eta_t&=& 0\eq
The above assumptions were chosen to insure these properties.\wwtbp
\par
In particular, the schemes given in (\ref{ab}) satisfy the hypotheses of the previous
proposition (notice that $b(U_\delta)\to b(U)$ as $\delta\to0$, due to the uniform convergence of $U_\delta$ to $U$), so that under the conditions of Theorem \ref{t1}, we get 
\bq
\lim_{t\ri +\iy} I_t&=&0\eq
Let us deduce (\ref{mN1})
for any neighborhood $\cN$ of the set $\cM$ of the global minima of $U$. From Cauchy-Schwartz inequality we have for any $t>0$,
\bq
\lVe m_t-\mu_{\beta_t,\delta_t}\rVe_{\mathrm{tv}}&=& \int \lve f_t-1\rve \, \mu_{\beta_t,\delta_t}\\
&\leq & \sqrt{I_t}\eq
An equivalent definition of the total variation norm states that
\bq
\lVe m_t-\mu_{\beta_t}\rVe_{\mathrm{tv}}&=&2\max_{A\in\cT}\lve m_t(A)-\mu_{\beta_t,\delta_t}(A)\rve\eq
where $\cT$ is the Borelian $\sigma$-algebra  of $M$.
It follows that (\ref{mN1}) reduces to \bq
\lim_{\beta,\delta^{-1}\ri +\iy} \mu_{\beta,\delta}(\cN)&=&1\eq
for any neighborhood $\cN$ of $\cM$ and $\delta$ sufficiently small, property which is immediate from the definition (\ref{mub}) of the Gibbs measures $\mu_{\beta,\delta}$ for $\beta\geq 0$ and $\delta>0$.
\begin{rem}
Similarly to the approach presented for instance in \cite{MR1275365,MR1425361},
we could have studied the evolution of $(E_t)_{t>0}$,  which are the
relative entropies of the time marginal laws with respect to the 
corresponding instantaneous Gibbs measures, namely
\bq
\fo t>0,\qquad E_t&\df& \int \ln\lt(\frac{m_t}{\mu_{\beta_t,\delta_t}}\rt)\, dm_t\eq
To get a differential 
  inequality satisfied by these functionals, the spectral gap estimate of 
\mbox{Holley}, Kusuoka and Stroock \cite{MR995752} recalled in Proposition \ref{HKS2}
must be replaced by the corresponding logarithmic Sobolev constant estimate.
\end{rem}

 \bibliographystyle{plain}

\begin{thebibliography}{}

\end{thebibliography}


\begin{thebibliography}{1}

\bibitem{Afsari-Tron-Vidal:11}
B. Afsari, R. Tron and R. Vidal, \emph{On the convergence of gradient descent for finding the Riemannian center of mass} arXiv:1201.0925

\bibitem{Arnaudon-Miclo:13}
Marc~Arnaudon and~Laurent~Miclo.
\newblock Means in complete manifolds: uniqueness and approximation.
\newblock Available on
  \texttt{http://hal.archives-ouvertes.fr/hal-00717677/}. To appear in ESAIM P.S.
  
\bibitem{Arnaudon}
Marc~Arnaudon and~Laurent~Miclo.
\newblock A stochastic algorithm finding $p$-means on the circle.
\newblock Available on
  \texttt{http://hal.archives-ouvertes.fr/hal-00781715}, 2013-01.
  
  \bibitem{Arnaudon-Nielsen:12a}
 M. Arnaudon, F. Nielsen, \emph{Medians and means in Finsler geometry}, LMS Journal of Computation and Mathematics volume 15, pp. 23-37 (2012)

\bibitem{Arnaudon-al:12}
{A}rnaudon, M., {D}ombry, C., {P}han, A., {Y}ang, L.,  \emph{{S}tochastic
  algorithms for computing means of probability measures}
 Stoch. Proc. Appl. 122 (2012), pp. 1437-1455.
 
 \bibitem{Arnaudon-Nielsen:12}
 M. Arnaudon, F. Nielsen, \emph{On computing the Riemannian 1-Center}, Computational Geometry: Theory and Applications, 46 (2013), no. 1, 93--104. 

\bibitem{Badoiu-Clarkson:03}
B\u{a}doiu, M., Clarkson, K.~L., 2003, \emph{Smaller core-sets for balls,} Proceedings of the fourteenth annual ACM-SIAM symposium on Discrete
  algorithms. Society for Industrial and Applied Mathematics, Philadelphia, PA,
  USA, pp. 801--802.

\bibitem{BP03}
R. Bhattacharya and V. Patrangenaru, \emph{Large sample theory of intrinsic and extrinsic sample means on manifolds (i)}, Annals of Statistics 31 (1), pp. 1--29 (2003)

\bibitem{Bonnabel:11}
S. Bonnabel, \emph{Convergence des m\'ethodes de gradient stochastique sur les variétés riemanniennes} In GRETSI. Bordeaux. 2011

\bibitem{Cardot-Cenac-Zitt:12}
H. Cardot, P. Cénac, P.-A. Zitt, \emph{Efficient and fast estimation of the geometric median in Hilbert spaces with an averaged stochastic gradient algorithm}, Bernoulli.

\bibitem{Charlier:11}
B. Charlier, \emph{Necessary and sufficient condition for the existence of a Fr\'echet mean on the circle} arXiv:1109.1986.

\bibitem{MR838085}
Stewart~N. Ethier and Thomas~G. Kurtz.
\newblock {\em Markov processes}.
\newblock Wiley Series in Probability and Mathematical Statistics: Probability
  and Mathematical Statistics. John Wiley \& Sons Inc., New York, 1986.
\newblock Characterization and convergence.

\bibitem{Groisser:05}
D. Groisser, \emph{Newton's method, zeroes of vector fields, and the Riemannian center of mass}
Adv. in Appl. Math. 33 (2004), no. 1, 95-135.

\bibitem{Groisser:06}
D. Groisser, \emph{On the convergence of some Procrustean averaging algorithms} Stochastics 77 (2005), no. 1, 31-60.

\bibitem{MR995752}
Richard~A. Holley, Shigeo Kusuoka, and Daniel~W. Stroock.
\newblock Asymptotics of the spectral gap with applications to the theory of
  simulated annealing.
\newblock {\em J. Funct. Anal.}, 83(2):333--347, 1989.

\bibitem{2011arXiv1108.2141H}
T.~{Hotz} and S.~{Huckemann}.
\newblock {Intrinsic Means on the Circle: Uniqueness, Locus and Asymptotics}.
\newblock {\em ArXiv e-prints}, August 2011.

\bibitem{Hsu:99}
E.P. Hsu
\newblock {Estimates of derivatives of the heat kernel on a compact Riemannian manifold}
\newblock {\em Proceedings of the American Mathematical Society.} Vol. 127 (1999), n. 12, pp 3739--3744.

\bibitem{MR602391}
Chii-Ruey Hwang.
\newblock Laplace's method revisited: weak convergence of probability measures.
\newblock {\em Ann. Probab.}, 8(6):1177--1182, 1980.

\bibitem{MR1011252}
Nobuyuki Ikeda and Shinzo Watanabe.
\newblock {\em Stochastic differential equations and diffusion processes},
  volume~24 of {\em North-Holland Mathematical Library}.
\newblock North-Holland Publishing Co., Amsterdam, second edition, 1989.

\bibitem{Le:04}
H. Le, \emph{Estimation of Riemannian barycentres}, LMS J. Comput. Math. 7 (2004), pp. 193--200

\bibitem{MR1275365}
Laurent Miclo.
\newblock Recuit simul\'e sans potentiel sur une vari\'et\'e riemannienne
  compacte.
\newblock {\em Stochastics Stochastics Rep.}, 41(1-2):23--56, 1992.

\bibitem{MR1348382}
Laurent Miclo.
\newblock Remarques sur l'ergodicit\'e des algorithmes de recuit simul\'e sur
  un graphe.
\newblock {\em Stochastic Process. Appl.}, 58(2):329--360, 1995.

\bibitem{MR1425361}
Laurent Miclo.
\newblock Recuit simul\'e partiel.
\newblock {\em Stochastic Process. Appl.}, 65(2):281--298, 1996.

\bibitem{MR2254442}
Xavier Pennec.
\newblock Intrinsic statistics on {R}iemannian manifolds: basic tools for
  geometric measurements.
\newblock {\em J. Math. Imaging Vision}, 25(1):127--154, 2006.

\bibitem{Sturm:05}
K.T. Sturm, \emph{Probability measures on metric spaces of nonpositive curvature}, Heat kernels and analysis on manifolds, graphs, and metric spaces (Paris, 2002), 357-390, Contemp. Math., 338, Amer. Math. Soc., Providence, RI, 2003.

\bibitem{Yang:10}
L.Yang, \emph{Riemannian median and its estimation}, 
LMS Journal of Computation and Mathematics, Vol 13 (2010), pp 461--479.


\end{thebibliography}
 
\def\cprime{$'$}

\vskip2cm
\hskip70mm\box5

\end{document}